\documentclass[11pt,twoside]{article}
\usepackage{amssymb,amsmath}
\usepackage{graphicx}
\usepackage{amsthm}
\newtheorem{theorem}{Theorem}
\newtheorem{lemma}{Lemma}

\newtheorem{definition}{Definition}
\newcommand{\zm}{{\bf Proof.}}
\newcommand{\zb}{\hfill$\blacksquare$}

\begin{document}

\thispagestyle{empty}\setcounter{page}{1}

\begin{center}
{\large A Weak Galerkin Finite Element Method for Burgers' Equation}\\[0.3cm]
\footnotetext{{\em Correspondence to}: Yanli Chen, Department of
Mathematics, Northeastern University,
Shenyang, 110004, China (e-mail: chenyanli@mail.neu.edu.cn)\\
Contract grant sponsor: Fundamental Research Funds for the Central Universities, No. 02060022116002; and the State Key Laboratory of Synthetical Automation for Process Industries Fundamental Research Funds, No. 2013ZCX02; and the National Natural Science Funds of China, No. 11371081.}
{Yanli Chen and Tie Zhang\\[0.2cm]}
{\em\small Department of
Mathematics, Northeastern University,
Shenyang, 110004, China}\\[0.4cm]
\end{center}
We propose a weak Galerkin(WG) finite element method for solving the one-dimensional Burgers' equation. Based on a new weak variational form, both semi-discrete and fully-discrete WG finite element schemes are established and analyzed. We prove the existence of the discrete solution and derive the optimal order error estimates in the discrete $H^1$-norm and $L^2$-norm, respectively. Numerical experiments are presented to illustrate our theoretical analysis. \\[0.3cm]

\section{Introduction}
\setcounter{section}{1}\setcounter{equation}{0}
In this paper, we consider the following Burgers' equation
\begin{subequations}\label{eq1}
\begin{equation}
\frac{\partial u}{\partial t}+u\frac{\partial u}{\partial x}-\nu\frac{\partial^2 u}{\partial x^2}=0,\quad x\in I,\quad t\in J,
\end{equation}
with the initial-boundary conditions
\begin{equation}
u(x,0)=g(x),\quad \forall x\in I,
\end{equation}
\begin{equation}
u(0,t)=0,\quad u(1,t)=0,\quad \forall t\in J,
\end{equation}
\end{subequations}
where $I=(0,1)$, $J=(0,T]$ with $T<\infty$, $u$ is the unknown velocity,
$\nu=\frac{1}{Re}$, ${ Re}$ is the Reynolds number and $g(x)$ is a given function. \par
Burgers' equation (\ref{eq1}a) has been widely
used for various applications, such as modeling of gas dynamics and turbulence \cite{b1}, describing wave propagation in acoustics
and hydrodynamics \cite{b2}, etc. More importantly, Burgers' equation is an appropriate simplified form of the
Navier-Stokes equations. Thus, it is anticipated that the crucial point in approximating the solution of the Navier-Stokes equations may be found
by numerical simulation of Burgers' equation. Furthermore, Burgers' equation will become a hyperbolic equation when the viscosity $\nu$ tends to zero. Therefore, a high performance numerical method for Burgers' equation may pave the way for solving first order hyperbolic equation. \par
Due to the importance of the Burgers' equation, various numerical methods were proposed
in the past decades. These methods mainly include finite difference, finite volume, finite element, boundary element and spectral methods, etc., see \cite{b3,b4,b28,b5,b6,b7,b8,b9,b10,b11,b12}, and the references therein.\par
Recently, the WG finite element method has attracted much attention in the field of numerical partial differential equations  \cite{b13,b14,b15,b16,b17,b18,b19,b20,b21,b22,b23,b24,b25}. This method was presented originally by Wang and Ye for solving general elliptic
problems in multi-dimensional domain \cite{b13}. Since then, some WG finite element methods have been developed to solve parabolic equations \cite{b15}, Stokes equations \cite{b20}and Biharmonic Equations \cite{b18,b25}, etc. In the present work, the WG finite element method is applied to one-dimensional Burgers' equation. In order to make the weak finite element equation have a positive property, we present a special weak form for problem (\ref{eq1}). Based on this weak form, we establish both semi-discrete and fully-discrete WG finite element schemes and derive the optimal order error estimates in the discrete $H^1$-norm, the $L^2$-norm, respectively.\par
An outline of the paper is as follows. In Section 2, we establish a special weak form for problem (\ref{eq1}) and introduce the definition of discrete weak derivative. Section 3 is devoted to semi-discrete WG finite element method. In Section 4, We construct fully-discrete WG finite element scheme and prove the optimal order error estimates. In Section 5, numerical experiments are presented to show the efficacy of the WG finite element method and confirm our theoretical analysis.\par
Throughout this paper, we adopt the notations $H^s(\tilde{I})$ to indicate the usual Sobolev spaces on subinterval $\tilde{I}\subset I$ equipped with the inner products $(\cdot,\cdot)_{s,\tilde{I}}$, the norms $\|\cdot\|_{s,\tilde{I}}$ and the seminorms $|\cdot|_{s,\tilde{I}}$, for any $s\geq0$. Denote by $(\cdot,\cdot)_{\tilde{I}}$ and $\|\cdot\|_{\tilde{I}}$ the $L^2(\tilde{I})$-inner product and $L^2(\tilde{I})$-norms, respectively. If $\tilde{I}=I$, we shall drop the subscript $\tilde{I}$ in the inner product and norm notations.\par
The space $H^l(0,T;H^s(I))$ is defined as
\begin{equation*}
    H^l(0,T;H^s(I))=\{v\in H^s(I); \int_0^T\sum\limits_{1\leq i\leq l}\Big{\|}\frac{d^iv(t)}{dt^i}\Big{\|}_{s,I}^2dt<\infty\},
\end{equation*}
equipped with the norm
\begin{equation*}
    \|v\|_{H^l(H^s(I))}=\Big{(}\int_0^T \sum\limits_{1\leq i\leq l}\Big{\|}\frac{d^iv(t)}{dt^i}\Big{\|}_{s,I}^2dt\Big{)}^{1/2}.
\end{equation*}
Especially, if $l = 0$, let
\begin{equation*}
    \|v\|_{L^2(H^s(I))}=\Big{(}\int_0^T\|v(t)\|_{s,I}^2dt\Big{)}^{1/2}.
\end{equation*}

\section{Weak Variational Form and Discrete Weak Derivative}
\setcounter{section}{2}\setcounter{equation}{0}
In order to define the weak Galerkin finite element approximation to problem (\ref{eq1}), we first need to derive a special weak variational form associated with problem (\ref{eq1}).\par
Multiplying equation (\ref{eq1}a) by $v\in H^1_0(I)$ and integrating, we obtain
\begin{equation}\label{eq11}
(u_t, v)+(uu_x,v)+\nu(u_x,v_x)=0.
\end{equation}
Integrating by parts, we have
\begin{eqnarray*}
   (uu_x,v)&=&\frac{1}{3}(uu_x,v)+\frac{1}{3}((u^2)_x,v)  \\
   &=&\frac{1}{3}(uu_x,v)-\frac{1}{3}(uu,v_x).
\end{eqnarray*}
Substituting this into equation (\ref{eq11}), we see that a weak form for problem (\ref{eq1}) is to find $u\in H^1(0,T; H^1_0(I))$ such that
\begin{equation}\label{eq13}
\left\{
\begin{aligned}
&(u_t, v)+\nu(u_x,v_x)+\frac{1}{3}(uu_x,v)-\frac{1}{3}(uu,v_x)=0,~\forall v\in H^1_0(I),~t\in J,\\
&u(x,0)=g(x),\quad x\in I.
\end{aligned}
\right.
\end{equation}
Let
\begin{equation*}
  a(w;u,v)=\nu(u_x,v_x)+\frac{1}{3}(wu_x,v)-\frac{1}{3}(wu,v_x).
\end{equation*}
Obviously, for any given $w$, $a(w;u,v)$ is positive define, i.e.,
\begin{equation*}
  a(w;u,u)\geq\nu \|u\|_1^2,\quad \forall u\in H^1_0(I).
\end{equation*}
\par
We now introduce the concepts of discrete weak derivative. Let closed interval $\bar{I}_a=[x_a,x_b]$ and its interior $I_a=(x_a,x_b)$. A weak function on $\bar{I}_a$ refers to a function $v=\{v^0,v^a,v^b\}$, where $v^0=v|_{I_a}$,$v^a=v(x_a)$ and $v^b=v(x_b)$. Note that $v^a$ and $v^b$ may not be necessarily the trace of $v^0$ at the interval endpoints $x_a$ and $x_b$. Define the weak function space by
\begin{equation*}
  W(I_a)=\{v=\{v^0,v^a,v^b\}:v^0\in L^2(I_a),~|v^a|+|v^b|<\infty\}.
\end{equation*}
\par
 For any integer $r\geq 0$, let $P_r(I_a)$ be the space composed of all polynomials on $I_a$ with degree no more than $r$.
\begin{definition}
For $v\in W(I_a)$, the discrete weak derivative $d_{w,r} v \in P_r(I_a)$ is defined as the unique solution of the following equation
\begin{equation}\label{eq6}
\int_{I_a}d_{w,r}vqdx=-\int_{I_a}v^0q'dx+v^bq^b-v^aq^a,\quad \forall q\in P_r(I_a),
\end{equation}
where $q^a=q(x_a)$, $q^b=q(x_b)$.
\end{definition}
\section{Semi-discrete WG Finite Element Method}
The goal of this section is to develop semi-discrete WG finite element method for problem (\ref{eq1}) and derive error estimates in the discrete $H^1$-norm and $L^2$-norm, respectively.\\
\subsection{Semi-discrete WG Finite Element Scheme}
Let $I_h:0=x_1<x_2\cdots<x_{N-1}<x_N=1$ be a partition of interval $I=(0,1)$, $I_i=(x_i,x_{i+1})$ be the partition element. Denote the meshsize by $h=\max\limits_i h_i$, where $h_i=x_{i+1}-x_i$, $i=1,2\cdots,N-1$. \par
For any given integer $k\geq 0$, the discrete weak function space associated with partition $I_h$ is defined by
\begin{equation}\label{eq7}
W(I_h,k)=\{v=\{v^0,v^L,v^R\}:v|_{I_i}\in W(I_i,k),~i=1,2,\cdots,N-1\},
\end{equation}
where
\begin{equation}\label{eq8}
W(I_i,k)=\{v=\{v^0,v^i,v^{i+1}\}:v^0\in P_k(I_i),~|v^i|+|v^{i+1}|<\infty\}.
\end{equation}
Note that for a weak function $v|_{I_i}\in W(I_i,k)$, the endpoint values $v^i$ and $v^{i+1}$ may be independent of the interior value $v^0$. Let $I_L=(x_{i-1},x_i)$ and $I_R=(x_i,x_{i+1})$ be two elements sharing the endpoints $x_i$. For a weak function $v\in W(I_h,k)$, set $v|_{\bar{I}_L}=\{v^0_L,v^{i-1}_L,v^i_L\}$ and $v|_{\bar{I}_R}=\{v^0_R,v^i_R,v^{i+1}_R\}$. We define the jump of weak function $v$ at point $x_i$ by
\begin{equation*}
    [v]_{x_i}=v^i_R-v^i_L, \quad v\in W(I_h,k).
\end{equation*}
It is obvious that $v$ is single value at at point $x_i$ if and only if $[v]_{x_i}=0$.
The weak finite element space $S_h$ is defined by
\begin{equation}\label{eq9}
S_h=\{v:v\in W(I_h,k),[v]_{x_i}=0,~i=1,2,\cdots,N-1\}.
\end{equation}
We refer to $S_h^0$ as the subspace of $S_h$ with vanish boundary value on endpoints of the interval $[0,1]$; i.e.,
\begin{equation}\label{eq10}
S^0_h=\{v:v\in S_h, v^1=0,v^N=0\}.
\end{equation}
Denote the discrete $L^2$ inner product and norm by
\begin{equation*}
  (u,v)_h=\sum_{i=1}^{N-1}(u,v)_{I_i}=\sum_{i=1}^{N-1}\int_{I_i}uvdx,\quad \|u\|_h^2=(u,u)_h.
\end{equation*}
\par
Based on the variational formulation (\ref{eq13}), we define semi-discrete WG finite element method by finding $u_h(t)\in S_h^0$ such that
\begin{equation}\label{eq14}
\left\{
\begin{aligned}
&\big{(}(u_h^0)_t, v^0\big{)}_h+\nu(d_{w,r}u_h,d_{w,r}v)_h+\frac{1}{3}(u_h^0d_{w,r}u_h,v^0)_h\\
&-\frac{1}{3}(u_h^0u_h^0,d_{w,r}v)_h=0,\quad\forall v\in S_h^0,~t\in J,\\
&u_h(x,0)=g_h(x),\quad x\in I,
\end{aligned}
\right.
\end{equation}
where $g_h\in S_h^0$ is a proper approximation of function $g$.\par
Taking $v=u_h$ in equation (\ref{eq14}), we get
\begin{equation*}
    \frac{1}{2}\frac{d}{dt}\|u_h^0\|_h^2+\nu\|d_{w,r}u_h\|_h^2=0.
\end{equation*}
Integrating with respect to t leads to
\begin{equation*}
    \|u_h^0(t)\|_h^2+2\nu\int_0^t\|d_{w,r}u_h(s)\|_h^2ds=\|u_h^0(0)\|_h^2,
\end{equation*}
which together with the continuation theorem guarantee the existence of solution $u_h$ of problem (\ref{eq14}).\\
\subsection{Error Analysis}
In this subsection, we shall derive error estimates for the semi-discrete WG finite element method in the discrete $H^1$-norm and $L^2$-norm, respectively.\par
To balance the approximation accuracy between spaces $S_h$ and
$P_r(I_i)$ used to compute $d_{w,r}u_h$, from now on, we always set the index $r = k+1$ in equations (\ref{eq6}) and (\ref{eq9}).\par
For $l\geq 0$, let $P_h^l$ is the local $L^2$ projection operator, restricted on each element $I_i$, $P_h^l: u\in L^2(I_i)\rightarrow P_h^lu\in P_l(I_i)$ such that
\begin{equation*}\label{eq21}
(u-P_h^lu,q)_{I_i}=0,\quad \forall q\in P_l(I_i),\quad i=1,2,\cdots,N-1.
\end{equation*}
By the Bramble-Hilbert lemma, it is easy to prove that
\begin{equation}\label{eq22}
\|u-P_h^lu\|_{I_i}+h_i\|u-P_h^lu\|_{1,I_i}\leq Ch_i^s\|u\|_{s,I_i}, \quad 0\leq s\leq l+1.
\end{equation}
We now define a projection operator $Q_h:u\in H^1(I)\rightarrow Q_hu \in W(I_h,k)$ such that
\begin{equation*}\label{eq23}
Q_hu|_{\bar{I}_i}=\{Q_h^0u,(Q_h u)^i,(Q_h u)^{i+1}\}=\{P_h^ku,u(x_i),u(x_{i+1})\},\quad i=1,2,\cdots,N-1.
\end{equation*}
Obviously, $Q_hu\in S_h^0$ if $u\in H^1_0(I)$. It follows from (\ref{eq22}) that
\begin{equation}\label{eq24}
\|Q_h^0 u-u\|_{I_i}\leq Ch_i^s\|u\|_{s,I_i},\quad 0\leq s\leq k+1.
\end{equation}
Furthermore, using the definition  of operator $Q_h$ and the discrete weak derivative $d_{w,r}$ in (\ref{eq6}), we have
\begin{eqnarray*}
  \int_{I_i}d_{w,r}Q_hu q dx &=& -\int_{I_i}Q_h^0uq^{\prime}dx+(Q_h u)^{i+1}q^{i+1}-(Q_h u)^iq^i\\
 &=&-\int_{I_i}uq^{\prime}dx+u^{i+1}q^{i+1}-u^iq^i\\
 &=&\int_{I_i}u^{\prime}qdx,\qquad\forall q\in P_r(I_i).
\end{eqnarray*}
Hence $d_{w,r}Q_h u=P_h^r u^{\prime}$ and (noting that $r=k+1$)
\begin{equation}\label{eq25}
\|d_{w,r}Q_h u-u^{\prime}\|_{I_i}\leq Ch_i^s\|u\|_{s+1,I_i},\quad0\leq s\leq k+2.
\end{equation}
Estimates (\ref{eq24}) and (\ref{eq25}) show that $Q_hu\in S_h^0$ is a very good approximation for function $u\in H_0^1(I)\cap H^{s+1}(I)$, $s\geq0$.\par
For error analysis, we still need to introduce the following projection function, which can be found in \cite{b26}.
\begin{lemma}\label{L:L1}
For $u\in H^1(I)$, there exists a projection function $\pi_h u\in H^1(I)$, restricted on element $I_i$, $\pi_h u\in P_r(I_i)$ satisfies(see \cite{b26} for details)
\begin{align}
\label{eq26}&\big{(}(\pi_hu)^{\prime},q\big{)}_{I_i}=(u^{\prime},q)_{I_i},\quad \forall q\in P_k(I_i),~i=1,2,\cdots,N-1, \\
 \label{eq27}            & \pi_hu(x_i)=u(x_i),\quad i=1,2,\cdots,N,\\
 \label{eq28}            &\|u-\pi_h u\|_{I_i}+h_i\|u^{\prime}-(\pi_hu)^{\prime}\|_{I_i}\leq C h_i^{s+1}\|u\|_{s+1},\quad 0\leq s \leq k+1.
\end{align}
\end{lemma}
\begin{lemma}\label{L:L2}
Let $u(t)\in H^1(0,T;H^2(I))$ be the solution of problem (\ref{eq1}). Then we have
\begin{equation}\label{eq20}
\begin{split}
&(u_t,v^0)_h+\nu(\pi_hu_x,d_{w,r}v)_h+\frac{1}{3}(uu_x,v^0)_h\\
&-\frac{1}{3}(\pi_hu^2,d_{w,r}v)_h=0, ~\forall v\in S_h^0,~t\in J.
\end{split}
\end{equation}
\end{lemma}\noindent
\zm\quad Multiplying the first equation of the system (\ref{eq1}) by $v^0$ and integrating, we obtain
\begin{equation}\label{eq2}
(u_t,v^0)_h+(uu_x,v^0)_h-\nu(u_{xx},v^0)_h=0, \quad\forall v\in S_h^0,~~t\in J.
\end{equation}
It follows from (\ref{eq26}) that
\begin{eqnarray}
\nonumber (uu_x,v^0)_h&=&\frac{1}{3}(uu_x,v^0)_h+\frac{1}{3}\big{(}(u^2)_x,v^0\big{)}_h\\\label{eq16}
&=&\frac{1}{3}(uu_x,v^0)_h+\frac{1}{3}\big{(}(\pi_h u^2)_x,v^0\big{)}_h
\end{eqnarray}
By the definition of operator $d_{w,r}$, we have
\begin{equation*}
    \big{(}(\pi_h u^2)_x,v^0\big{)}_{I_i}=-(\pi_h u^2,d_{w,r}v)_{I_i}+(\pi_h u^2)^{i+1}v^{i+1}-(\pi_h u^2)^iv^i,\quad\forall v\in S_h^0,
\end{equation*}
Summing and noting that $v^1=0$ and $v^N=0$, it leads to
\begin{equation}\label{eq3}
    \big{(}(\pi_h u^2)_x,v^0\big{)}_h=-(\pi_h u^2,d_{w,r}v)_h.
\end{equation}
Combining (\ref{eq16}) with (\ref{eq3}) yields
\begin{equation}\label{eq4}
(uu_x,v^0)_h=\frac{1}{3}(uu_x,v^0)_h-\frac{1}{3}(\pi_h u^2,d_{w,r}v)_h
\end{equation}
Similarly, following equation \eqref{eq26} and the definition of operator $d_{w,r}$, for any $v\in S_h^0$, it holds
\begin{equation}\label{eq5}
-(u_{xx},v^0)_h=-\big{(}(\pi_h u_x)_x,v^0\big{)}_h=(\pi_h u_x, d_{w,r} v)_h,
\end{equation}
Substituting (\ref{eq4}) and (\ref{eq5}) into (\ref{eq2}), we complete the proof.\zb

\begin{theorem}\label{T:T1}
Let $u(x,t)$ and $u_h(x,t)$ be the solution of problems (\ref{eq1}) and (\ref{eq14}), respectively, $u\in H^1(0,T; H^{1+r}(I))$ and $g_h=Q_hg$. Then, there exists a constant $C$ such that
\begin{eqnarray*}
   &&\|u-u_h^0\|_h^2+\nu\int_0^t\|u_x(s)-d_{w,r}u_h(s)\|_h^2ds \\
 &\leq& Ch^{2(k+1)}\Big{(}\|u\|_{k+1}^2+\int_0^t\|u(s)\|^2_{k+2}ds\Big{)}.
\end{eqnarray*}
\end{theorem}
\noindent\zm\quad Subtracting (\ref{eq14}) from (\ref{eq20}), and using $(Q_h^0u_t,v^0)_h=(u_t,v^0)_h$, we arrive at
\begin{eqnarray}
\nonumber&&\big{(}Q_h^0u_t-(u_h^0)_t,v^0\big{)}_h+\nu\big{(}d_{w,r}(Q_hu-u_h),d_{w,r}v\big{)}_h \\ \nonumber&=&\frac{1}{3}(u_h^0d_{w,r}u_h,v^0)_h-\frac{1}{3}(uu_x,v^0)_h\\
\nonumber &+&\frac{1}{3}(\pi_hu^2,d_{w,r}v)_h-\frac{1}{3}(u_h^0u_h^0,d_{w,r}v)_h\\\label{eq18}
&+&\nu(d_{w,r}Q_hu,d_{w,r}v)_h-\nu(\pi_hu_x,d_{w,r}v)_h
\end{eqnarray}
Let $e_h=Q_hu-u_h$. Taking $v=e_h$ in equation (\ref{eq18}), we have
\begin{eqnarray}
\nonumber&&\frac{1}{2}\frac{d}{dt}\|e_h^0\|_h^2+\nu\|d_{w,r}e_h\|_h^2 \\
\nonumber&=&\frac{1}{3}(u_h^0d_{w,r}u_h,e_h^0)_h-\frac{1}{3}(uu_x,e_h^0)_h\\
\nonumber &+&\frac{1}{3}(\pi_hu^2,d_{w,r}e_h)_h-\frac{1}{3}(u_h^0u_h^0,d_{w,r}e_h)_h\\
\nonumber&+&\nu(d_{w,r}Q_hu,d_{w,r}e_h)_h-\nu(\pi_hu_x,d_{w,r}e_h)_h\\\label{eq17}
&=&E_1+E_2+E_3.
\end{eqnarray}
In what follows, we estimate $E_1$, $E_2$ and $E_3$ in equation \eqref{eq17}. \par
To estimate $E_1$ and $E_2$, we can rewritten $E_1$ and $E_2$ as follows:
\begin{eqnarray*}
 E_1&=&\frac{1}{3}\Big{(}\big{(}u_h^0d_{w,r}(u_h-Q_hu),e_h^0\big{)}_h+\big{(}(u_h^0-Q_h^0u)d_{w,r}Q_hu,e_h^0\big{)}_h\\
 &&+\big{(}(Q_h^0u-u)d_{w,r}Q_hu,e_h^0\big{)}_h+\big{(}u(d_{w,r}Q_hu-u_x),e_h^0\big{)}_h\Big{)},
\end{eqnarray*}
\begin{eqnarray*}
E_2 &=&\frac{1}{3}\Big{(}(\pi_hu^2-u^2,d_{w,r}e_h)_h+\big{(}(u+Q_h^0u)(u-Q_h^0u),d_{w,r}e_h\big{)}_h\\
  &&+\big{(}Q_h^0u(Q_h^0u-u_h^0),d_{w,r}e_h\big{)}_h+\big{(}u_h^0(Q_h^0u-u_h^0),d_{w,r}e_h\big{)}_h\Big{)}
\end{eqnarray*}
Note that the first term in $E_1$ is the same as the last term in $E_2$ except for the sign. Thus, we can obtain that
\begin{eqnarray}
\nonumber&& E_1+E_2 \\
\nonumber &=&\frac{1}{3}\Big{(}(-e_h^0d_{w,r}Q_hu,e_h^0)_h+\big{(}(Q_h^0u-u)d_{w,r}Q_hu,e_h^0\big{)}_h\\
 \nonumber &+&\big{(}u(d_{w,r}Q_hu-u_x),e_h^0\big{)}_h+(\pi_hu^2-u^2,d_{w,r}e_h)_h\\\label{eq53}
  &+&\big{(}(u+Q_h^0u)(u-Q_h^0u),d_{w,r}e_h\big{)}_h+(Q_h^0ue_h^0,d_{w,r}e_h)_h \Big{)}.
\end{eqnarray}
It follows from $\varepsilon$-inequality and $d_{w,r}Q_hu=P_h^ru_x$ that
\begin{eqnarray}
\nonumber |E_1+E_2|&\leq&\frac{1}{3}\Big{(}\|u_x\|^2_{\infty}\|e_h^0\|_h^2+\frac{1}{2}\|u_x\|^2_{\infty}\|u-Q_h^0u\|_h^2+\frac{1}{2}\|e_h^0\|_h^2\\
\nonumber  &+&\frac{1}{2}\|u\|^2_{\infty}\|d_{w,r}Q_hu-u_x\|_h^2+\frac{1}{2}\|e_h^0\|_h^2\\
 \nonumber &+&\frac{1}{\nu}\|\pi_hu^2-u^2\|_h^2+\frac{\nu}{4}\|d_{w,r}e_h\|_h^2\\
 \nonumber &+&\frac{1}{\nu}\|u\|^2_{\infty}\|u-Q_h^0u\|_h^2+\frac{\nu}{4}\|d_{w,r}e_h\|_h^2\\\label{eq36}
  &+&\frac{1}{\nu}\|u\|^2_{\infty}\|e_h^0\|_h^2+\frac{\nu}{4}\|d_{w,r}e_h\|_h^2\Big{)}.
\end{eqnarray}
To estimate $E_3$, we rewritten $E_3$ as follow:
\begin{equation*}
    E_3=\nu(d_{w,r}Q_hu-u_x,d_{w,r}v)_h+\nu(u_x-\pi_h u_x,d_{w,r}v)_h.
\end{equation*}
Using $\varepsilon$-inequality, we have
\begin{eqnarray}
 \nonumber   |E_3|&\leq&2\nu\|d_{w,r}Q_hu-u_x\|_h^2+\frac{\nu}{8}\|d_{w,r}e_h\|_h^2 \\\label{eq37}
  &+&2\nu\|u_x-\pi_hu_x\|_h^2+\frac{\nu}{8}\|d_{w,r}e_h\|_h^2.
\end{eqnarray}
Substituting (\ref{eq36}) and (\ref{eq37}) into (\ref{eq17}), we get
\begin{eqnarray*}
 \nonumber&&\frac{1}{2}\frac{d}{dt}\|e_h^0\|_h^2+\frac{\nu}{2}\|d_{w,r}e_h\|_h^2\\
 \nonumber &\leq&\frac{1}{3}(\|u_x\|_{\infty}^2+1+\frac{1}{\nu}\|u\|^2_{\infty})\|e_h^0\|_h^2+\frac{1}{3}(\frac{1}{2}\|u_x\|^2_{\infty}+\frac{1}{\nu}\|u\|^2_{\infty})\|u-Q_h^0u\|_h^2\\
  \nonumber &+&(\frac{1}{6}\|u\|^2_{\infty}+2\nu)\|d_{w,r}Q_hu-u_x\|_h^2+\frac{1}{3\nu}\|\pi_hu^2-u^2\|_h^2+2\nu\|u_x-\pi_hu_x\|_h^2,
\end{eqnarray*}
together with (\ref{eq24}), (\ref{eq25}) and (\ref{eq28}), it yields
\begin{equation*}
    \frac{1}{2}\frac{d}{dt}\|e_h^0\|_h^2+\frac{\nu}{2}\|d_{w,r}e_h\|_h^2\leq C(\|e_h^0\|_h^2+h^{2(k+1)}\|u\|_{k+2}^2).
\end{equation*}
Integrating with respect to $t$ (noting that $e_h(0)=0$) and using the Gronwall lemma, we can get the following estimate
\begin{equation}\label{eq40}
\|e_h^0\|_h^2+\nu\int_0^t\|d_{w,r}e_h(s)\|_h^2ds\leq Ch^{2(k+1)}\int_0^t\|u(s)\|^2_{k+2}ds.
\end{equation}
Then, by the triangle inequality, approximation property \eqref{eq24} and \eqref{eq25}, we have
\begin{eqnarray*}
   &&\|u-u_h^0\|_h^2+\nu\int_0^t\|u_x(s)-d_{w,r}u_h(s)\|_h^2ds\\
   &\leq&C\Big{(}\|u-Q_h^0u\|_h^2+\|Q_h^0u-u_h^0\|_h^2+\nu\int_0^t\|u_x(s)-d_{w,r}Q_hu(s)\|_h^2ds\\
   &+&\nu\int_0^t\|d_{w,r}Q_hu(s)-d_{w,r}u_h(s)\|_h^2ds\Big{)}\\
   &\leq&C\Big{(}h^{2(k+1)}\|u\|^2_{k+1}+h^{2(k+1)}\int_0^t\|u(s)\|^2_{k+2}ds\Big{)},
\end{eqnarray*}
which completes the proof.\zb

\begin{theorem}\label{T:T2}
Let $u(x,t)$ and $u_h(x,t)$ be the solution of problems (\ref{eq1}) and (\ref{eq14}), respectively, $u\in H^1(0,T; H^{1+r}(I))$ and $g_h=Q_hg$. Then, there exists a constant $C$ such that
\begin{equation}\label{eq12}
  \|u_x-d_{w,r}u_h\|_h^2 \leq Ch^{2(k+1)}\Big{(}\|u(0)\|_{k+2}^2+\int_0^t\|u_t(s)\|^2_{k+2}ds\Big{)}.
\end{equation}
\end{theorem}
\noindent\zm\quad Taking $v=(e_h)_t$ in equation (\ref{eq18}), we have
\begin{eqnarray}
\nonumber&&\|(e_h^0)_t\|_h^2+\frac{\nu}{2}\frac{d}{dt}\|d_{w,r}e_h\|_h^2 \\
\nonumber&=&\frac{1}{3}\big{(}u_h^0d_{w,r}u_h,(e_h^0)_t\big{)}_h-\frac{1}{3}\big{(}uu_x,(e_h^0)_t\big{)}_h\\
\nonumber&+&\frac{1}{3}\big{(}\pi_hu^2,d_{w,r}(e_h)_t\big{)}_h-\frac{1}{3}\big{(}u_h^0u_h^0,d_{w,r}(e_h)_t\big{)}_h\\
\nonumber&+&\nu\big{(}d_{w,r}Q_hu,d_{w,r}(e_h)_t\big{)}_h-\nu\big{(}\pi_hu_x,d_{w,r}(e_h)_t\big{)}_h,\\\label{eq52}
&=&R_1+R_2+R_3.
\end{eqnarray}
On the analogy of (\ref{eq53}), we can obtain that
\begin{eqnarray}
\nonumber   &&R_1+R_2\\
\nonumber  &=&\frac{1}{3}\big{(}-e_h^0d_{w,r}Q_h u,(e_h^0)_t\big{)}_h+\frac{1}{3}\big{(}(Q_h^0u-u)d_{w,r}Q_hu,(e_h^0)_t\big{)}_h\\
\nonumber  &+&\frac{1}{3}\big{(}u(d_{w,r}Q_hu-u_x),(e_h^0)_t\big{)}_h+\frac{1}{3}\big{(}\pi_hu^2-u^2,d_{w,r}(e_h)_t\big{)}_h\\
\nonumber &+&\frac{1}{3}\big{(}(u+Q_h^0u)(u-Q_h^0u),d_{w,r}(e_h)_t\big{)}_h+\frac{1}{3}\big{(}Q_h^0ue_h^0,d_{w,r}(e_h)_t\big{)}_h\\\label{eq41}
&=&W_1+W_2+W_3+W_4+W_5+W_6
\end{eqnarray}
Using $\varepsilon$-inequality, (\ref{eq24}), (\ref{eq25}) and $d_{w,r}Q_hu=P_h^ru_x$, we can bound $W_1$, $W_2$ and $W_3$ as follows:
\begin{eqnarray}
\nonumber W_1&\leq&\frac{1}{3}\Big{(}\frac{1}{3}\|u_x\|_{\infty}^2\|e_h^0\|_h^2+\frac{3}{4}\|(e_h^0)_t\|_h^2\Big{)}\\\label{eq42}
  &\leq& \frac{1}{9}\|u_x\|_{\infty}^2\|e_h^0\|_h^2+\frac{1}{4}\|(e_h^0)_t\|_h^2
\end{eqnarray}
\begin{eqnarray}
\nonumber W_2 &\leq&\frac{1}{3}\Big{(}\frac{1}{3}\|u_x\|_{\infty}^2\|Q_h^0u-u\|_h^2+\frac{3}{4}\|(e_h^0)_t\|_h^2\Big{)} \\\label{eq43}
 &\leq&Ch^{2(k+1)}\|u\|_{k+1}^2+\frac{1}{4}\|(e_h^0)_t\|_h^2,
\end{eqnarray}
\begin{eqnarray}
 \nonumber W_3 &\leq&\frac{1}{3}\Big{(}\frac{1}{3}\|u\|_{\infty}^2\|d_{w,r}Q_hu-u_x\|_h^2+\frac{3}{4}\|(e_h^0)_t\|_h^2\Big{)} \\\label{eq44}
 &\leq& Ch^{2(k+1)}\|u\|_{k+2}^2+\frac{1}{4}\|(e_h^0)_t\|_h^2,
\end{eqnarray}
In the following, we estimate $W_4$, $W_5$ and $W_6$ in \eqref{eq41}. For the term $W_4$, we use Cauchy inequality and approximation property (\ref{eq28}) to obtain
\begin{eqnarray*}
\nonumber W_4&=&\frac{1}{3}\frac{d}{dt}(\pi_h u^2-u^2,d_{w,r}e_h)_h-\frac{1}{3}\big{(}\frac{d}{dt}(\pi_hu^2-u^2),d_{w,r}e_h\big{)}_h \\
  \nonumber &\leq&\frac{1}{3}\frac{d}{dt}(\pi_h u^2-u^2,d_{w,r}e_h)_h+\frac{1}{6}\|\frac{d}{dt}(\pi_hu^2-u^2)\|^2_h+\frac{1}{6}\|d_{w,r}e_h\|_h^2\\
   &\leq&\frac{1}{3}\frac{d}{dt}(\pi_h u^2-u^2,d_{w,r}e_h)_h+Ch^{2(k+1)}\|u_t\|_{k+1}^2+\frac{1}{6}\|d_{w,r}e_h\|_h^2.
\end{eqnarray*}
Then, the terms $W_5$ and $W_6$ can be bounded similarly as follows:
\begin{eqnarray*}
\nonumber W_5&\leq&\frac{1}{3}\frac{d}{dt}\big{(}(u+Q_h^0u)(u-Q_h^0u),d_{w,r}e_h\big{)}_h+Ch^{2(k+1)}\|u\|_{k+1}^2\\
&+&Ch^{2(k+1)}\|u_t\|_{k+1}^2+\frac{1}{3}\|d_{w,r}e_h\|_h^2,
\end{eqnarray*}
\begin{eqnarray*}
\nonumber W_6&\leq&\frac{1}{3}\frac{d}{dt}(Q_h^0ue_h^0,d_{w,r}e_h)_h+\frac{1}{6}\|u_t\|_{\infty}^2\|e_h^0\|_h^2+\frac{1}{6}\|d_{w,r}e_h\|_h^2\\
   &+&\frac{1}{4}\|(e_h^0)_t\|_h^2+\frac{1}{9}\|u\|_{\infty}^2\|d_{w,r}e_h\|_h^2.
\end{eqnarray*}
Substituting the estimates for $W_1 \sim W_6$ into (\ref{eq41}), we arrive at
\begin{eqnarray}
 \nonumber R_1+R_2 &\leq&(\frac{1}{9}\|u_x\|_{\infty}^2+\frac{1}{6}\|u_t\|_{\infty}^2)\|e_h^0\|_h^2+Ch^{2(k+1)}\|u\|_{k+2}^2 \\
 \nonumber &+&Ch^{2(k+1)}\|u_t\|_{k+1}^2+(\frac{2}{3}+\frac{1}{9}\|u\|_{\infty}^2)\|d_{w,r}e_h\|_h^2\\
 \nonumber &+&\|(e_h^0)_t\|_h^2 +\frac{1}{3}\frac{d}{dt}(\pi u^2-u^2,d_{w,r}e_h)_h\\
 \nonumber &+&\frac{1}{3}\frac{d}{dt}\big{(}(u+Q_h^0u)(u-Q_h^0u),d_{w,r}e_h\big{)}_h\\\label{eq54}
  &+&\frac{1}{3}\frac{d}{dt}(Q_h^0ue_h^0,d_{w,r}e_h)_h.
\end{eqnarray}
Similarly, we obtain the estimate for $R_3$,
\begin{eqnarray}
\nonumber   R_3&\leq& Ch^{2(k+1)}\|u_t\|_{k+2}^2+\nu\|d_{w,r}e_h\|_h^2 \\
\nonumber   &+&\nu\frac{d}{dt}(d_{w,r}Q_hu-u_x,d_{w,r}e_h)_h\\\label{eq55}
   &+&\nu\frac{d}{dt}(u_x-\pi_hu_x,d_{w,r}e_h)_h.
\end{eqnarray}
It follows from (\ref{eq52}), (\ref{eq54}) and (\ref{eq55}) that
\begin{eqnarray*}
  &&\frac{\nu}{2} \frac{d}{dt}\|d_{w,r}e_h\|_h^2\\
  &\leq&(\frac{1}{9}\|u_x\|_{\infty}^2+\frac{1}{6}\|u_t\|_{\infty}^2)\|e_h^0\|_h^2+Ch^{2(k+1)}\|u\|_{k+2}^2\\
  &+&Ch^{2(k+1)}\|u_t\|_{k+2}^2+(\frac{2}{3}+\frac{1}{9}\|u\|_{\infty}^2+\nu)\|d_{w,r}e_h\|_h^2\\
  &+&\frac{1}{3}\frac{d}{dt}(\pi u^2-u^2,d_{w,r}e_h)_h+\frac{1}{3}\frac{d}{dt}\big{(}(u+Q_h^0u)(u-Q_h^0u),d_{w,r}e_h\big{)}_h\\
  &+&\frac{1}{3}\frac{d}{dt}(Q_h^0ue_h^0,d_{w,r}e_h)_h+\nu\frac{d}{dt}(d_{w,r}Q_hu-u_x,d_{w,r}e_h)_h\\
&+&\nu\frac{d}{dt}(u_x-\pi_hu_x,d_{w,r}e_h)_h
\end{eqnarray*}
Thus, integrating with respect to $t$ (noting that $d_{w,r}e_h(0)=0$) and using estimate \eqref{eq40}, we have
\begin{eqnarray*}
\frac{\nu}{2} \|d_{w,r}e_h\|_h^2&\leq&C\Big{(}h^{2(k+1)}\int_0^t\|u(s)\|^2_{k+2}ds+h^{2(k+1)}\int_0^t\|u_t(s)\|^2_{k+2}ds\Big{)}\\
&+&\frac{1}{3}(\pi u^2-u^2,d_{w,r}e_h)_h+\frac{1}{3}\big{(}(u+Q_h^0u)(u-Q_h^0u),d_{w,r}e_h\big{)}_h\\
&+&\frac{1}{3}(Q_h^0ue_h^0,d_{w,r}e_h)_h+\nu(d_{w,r}Q_hu-u_x,d_{w,r}e_h)_h\\
&+&\nu(u_x-\pi_hu_x,d_{w,r}e_h)_h
   \end{eqnarray*}
Using $\varepsilon$-inequality, (\ref{eq24}), (\ref{eq25}) and (\ref{eq28}), it yields
\begin{eqnarray*}
\|d_{w,r}e_h\|_h^2  &\leq&C\Big{(}h^{2(k+1)}\|u\|_{k+2}^2+h^{2(k+1)}\int_0^t\|u(s)\|^2_{k+2}ds \\
&+&h^{2(k+1)}\int_0^t\|u_t(s)\|^2_{k+2}ds\Big{)}\\
&\leq&C\Big{(}h^{2(k+1)}\|u(0)\|_{k+2}^2+h^{2(k+1)}\int_0^t\|u_t(s)\|^2_{k+2}ds\Big{)}.
\end{eqnarray*}
Then, by the triangle inequality, we have the desired error estimate \eqref{eq12}.\zb
\section{Fully-discrete WG Finite Element Method}
In this section, we shall establish the fully-discrete WG Finite Element Method for problem \eqref{eq1} and derive error estimates in the discrete $H^1$-norm and $L^2$-norm, respectively.
\subsection{Fully-discrete WG Finite Element Scheme}
Let $0=t_0<t_1<\cdots<t_M=T$ be a partition for time interval $[0,T]$, where $t_n=n\tau, \tau=T/M$, $M$ is a positive integer. We refer to $U_n\in S_h^0$ as the approximation of $u(t_n)$ to be determined. Using a backward difference quotient $\bar{\partial}_tU_n=(U_n-U_{n-1})/\tau$ to approximate $(u_h)_t$ in problem \eqref{eq14}, we obtain the fully-discrete WG finite element scheme: find $U_n\in S_h^0~(n=1,2,\cdots,M)$ such that
\begin{equation}\label{eq15}
\left\{
\begin{aligned}
&(\bar{\partial}_tU_n^0, v^0)_h+\nu(d_{w,r}U_n,d_{w,r}v)_h+\frac{1}{3}(U_n^0d_{w,r}U_n,v^0)_h\\
&-\frac{1}{3}(U_n^0U_n^0,d_{w,r}v)_h=0,~\forall v\in S_h^0,~n=1,\cdots,M,\\
&U_0=g_h,\quad x\in I,
\end{aligned}
\right.
\end{equation}
where $g_h$ is a proper approximation of function $g$.
\begin{lemma}\label{L:L3}
There exists a solution $U_n\in S_h^0,\quad n=1,2,\cdots,M$ satisfying problem (\ref{eq15}).
\end{lemma}
\noindent\zm\quad Let $U_{n-1}$ be given. Define mapping $F: W_n\in S_h^0\rightarrow U_n=F(W_n)$ such that
\begin{equation}\label{eq29}
\left\{
\begin{aligned}
&(\bar{\partial}_tU_n^0, v^0)_h+\nu(d_{w,r}U_n,d_{w,r}v)_h+\frac{1}{3}(W_n^0d_{w,r}U_n,v^0)_h\\
&-\frac{1}{3}(W_n^0U_n^0,d_{w,r}v)_h=0,~\forall v\in S_h^0,~n\geq0\\
&U_0=g_h\quad x\in I.
\end{aligned}
\right.
\end{equation}
Taking $v=U_n$ in \eqref{eq29} and using Cauchy inequality, we have
\begin{equation}\label{eq19}
  \|U_n^0\|_h^2+2\tau\nu\|d_{w,r}U_n\|_h^2\leq\|U_{n-1}^0\|_h^2.
\end{equation}
This stability shows that for given $W_n\in S_h^0$, the linear system of equations \eqref{eq29} has a unique solution $U_n\in S_h^0$; that is, the mapping $F$ is well posed.\par
Let bounded set $D_h=\{U_n\in S_h^0:\|U_n^0\|_h^2+2\tau\nu\|d_{w,r}U_n\|_h^2\leq\|U_{n-1}^0\|_h^2\}$. The estimate \eqref{eq19} shows that $F$ is a mapping from $D_h$ into $D_h$. So according to the Brower fixed theorem,  for given $U_{n-1}$, problem \eqref{eq15} has a solution $U_n\in S_h^0$. \zb
\subsection{Error Analysis}
We are ready to derive some error estimates for the fully-discrete time WG approximation $U_n$ as the following theorems.
\begin{theorem}\label{T:T3}
Let $u$ and $U_n$ be the solution of problems (\ref{eq1}) and (\ref{eq15}), respectively, $u\in H^2(0,T;H^{1+r}(I))$ and $g_h=Q_hg$. Then, there exists a constant $C$ such that
\begin{equation*}
    \|u(t_n)-U_n^0\|_h^2\leq C\Big{(}h^{2(k+1)}\max\limits_{1\leq j\leq n}\|u(t_j)\|_{k+2}^2+{\tau}^2\int_0^{t_n}\|u_{tt}(s)\|^2ds\Big{)}.
\end{equation*}

\end{theorem}
\noindent\zm\quad Set $t=t_n$ in equation (\ref{eq20}) to obtain
\begin{eqnarray}
   \nonumber&& \big{(}u_t(t_n),v^0\big{)}_h+\nu\big{(}\pi_h u_x(t_n),d_{w,r}v\big{)}_h \\\label{q13}
  &&+\frac{1}{3}\big{(}u(t_n)u_x(t_n),v^0\big{)}_h-\frac{1}{3}\big{(}\pi_hu^2(t_n),d_{w,r}v\big{)}_h=0.
\end{eqnarray}
Subtracting (\ref{eq15}) from (\ref{q13}) and using $\big{(}\bar{\partial}_tQ_h^0u(t_n), v^0\big{)}_h=\big{(}\bar{\partial}_tu(t_n), v^0\big{)}_h$, we have
\begin{eqnarray}
\nonumber&& \big{(}\bar{\partial}_t(Q_h^0u(t_n)-U^0_n), v\big{)}_h+\nu\big{(}d_{w,r}(Q_hu(t_n)-U_n),d_{w,r}v\big{)}_h  \\
\nonumber  &=&\frac{1}{3}(U^0_nd_{w,r}U_n,v^0)_h-\frac{1}{3}\big{(}u(t_n)u_x(t_n),v^0\big{)}_h\\
 \nonumber &+&\frac{1}{3}\big{(}\pi_hu^2(t_n),d_{w,r}v\big{)}_h-\frac{1}{3}(U^0_nU^0_n,d_{w,r}v)_h\\
\nonumber  &+&\nu\big{(}d_{w,r}Q_hu(t_n),d_{w,r}v\big{)}_h-\nu\big{(}\pi_h u_x(t_n),d_{w,r}v\big{)}_h\\\label{q3}
      &+&\big{(}\bar{\partial}_tu(t_n)-u_t(t_n),v^0\big{)}_h
\end{eqnarray}
Let $e_n=Q_hu(t_n)-U_n$, and taking $v=e_n$ in (\ref{q3}), we arrive at
\begin{eqnarray}
\nonumber&& \frac{1}{\tau}\|e_n^0\|_h^2-\frac{1}{\tau}(e^0_{n-1}, e_n^0)_h+\nu\|d_{w,r}e_n\|_h^2  \\
\nonumber  &=&\frac{1}{3}(U_n^0d_{w,r}U_n,e_n^0)_h-\frac{1}{3}\big{(}u(t_n)u_x(t_n),e_n^0\big{)}_h\\
 \nonumber &+&\frac{1}{3}\big{(}\pi_hu^2(t_n),d_{w,r}e_n\big{)}_h-\frac{1}{3}(U_n^0U_n^0,d_{w,r}e_n)_h\\
\nonumber  &+&\nu\big{(}d_{w,r}Q_hu(t_n),d_{w,r}e_n\big{)}_h-\nu\big{(}\pi_h u_x(t_n),d_{w,r}e_n\big{)}_h\\
\nonumber      &+&\big{(}\bar{\partial}_tu(t_n)-u_t(t_n),e_n^0\big{)}_h\\\label{q1}
      &=&F_1+F_2+F_3+F_4.
\end{eqnarray}
Similar to the argument we make in Theorem \ref{T:T1}, the first three terms on the right hand side of (\ref{q1}) can be estimated as follows:
\begin{eqnarray}
  \nonumber|F_1+F_2|&\leq& \frac{1}{3}\Big{(}\|u_x(t_n)\|_{\infty}^2\|e_n^0\|_h^2+\frac{1}{2}\|u_x(t_n)\|_{\infty}^2\|Q_h^0u(t_n)-u(t_n)\|_h^2+\frac{1}{2}\|e_n^0\|_h^2\\
\nonumber &+&\frac{1}{2}\|u(t_n)\|_{\infty}^2\|d_{w,r}Q_hu(t_n)-u_x(t_n)\|_h^2+\frac{1}{2}\|e_n^0\|_h^2\\
\nonumber&+&\frac{1}{2\nu}\|\pi_h u^2(t_n)-u^2(t_n)\|_h^2+\frac{\nu}{2}\|d_{w,r}e_n\|_h^2\\
 \nonumber&+&\frac{1}{2\nu}\|u(t_n)\|_{\infty}^2\|u(t_n)-Q_h^0u(t_n)\|_h^2+\frac{\nu}{2}\|d_{w,r}e_n\|_h^2\\\label{q7}
&+&\frac{1}{2\nu}\|u(t_n)\|_{\infty}^2\|e_n^0\|_h^2+\frac{\nu}{2}\|d_{w,r}e_n\|_h^2\Big{)},
\end{eqnarray}
\begin{eqnarray}
 \nonumber  |F_3|&\leq&\nu\|d_{w,r}Q_hu(t_n)-u_x(t_n)\|_h^2+\frac{\nu}{4}\|d_{w,r}e_n\|_h^2  \\\label{q8}
   &+& \nu\|u_x(t_n)-\pi_h u_x(t_n)\|_h^2+\frac{\nu}{4}\|d_{w,r}e_n\|_h^2.
\end{eqnarray}
For the term $F_4$, we have
\begin{equation}\label{q9}
|F_4|\leq\frac{1}{2}\|\bar{\partial}_tu(t_n)-u_t(t_n)\|_h^2+\frac{1}{2}\|e_n^0\|_h^2.
\end{equation}
Substituting (\ref{q7}), (\ref{q8}) and (\ref{q9}) into (\ref{q1}), and using $\varepsilon$-inequality, we obtain
\begin{eqnarray}
 \nonumber  \|e_n^0\|_h^2 &\leq&\|e^0_{n-1}\|_h^2+2\tau\Big{(}\frac{1}{3}\|u_x(t_n)\|_{\infty}^2+\frac{5}{6}+\frac{1}{6\nu}\|u(t_n)\|_{\infty}^2\Big{)}\|e_n^0\|_h^2 \\
\nonumber    &&+2\tau\Big{(}\frac{1}{6}\|u_x(t_n)\|_{\infty}^2+\frac{1}{6\nu}\|u(t_n)\|_{\infty}^2\Big{)}\|Q_h^0u(t_n)-u(t_n)\|_h^2\\
\nonumber    &&+2\tau\Big{(}\frac{1}{6}\|u(t_n)\|_{\infty}^2+\nu\Big{)}\|d_{w,r}Q_hu(t_n)-u_x(t_n)\|_h^2\\
  \nonumber &&+\frac{\tau}{3\nu}\|\pi_h u^2(t_n)-u^2(t_n)\|_h^2+2\tau\nu\|u_x(t_n)-\pi_h u_x(t_n)\|_h^2\\\label{q10}
   &&+\tau\|\bar{\partial}_tu(t_n)-u_t(t_n)\|_h^2.
\end{eqnarray}
By repeated application(noting that $e_0=0$), and together with (\ref{eq24}), (\ref{eq25}), and (\ref{eq28}), we arrive at
\begin{eqnarray}
\nonumber  \|e_n^0\|_h^2 &\leq&C\tau \sum_{j=1}^n\|e^0_j\|_h^2+C\tau h^{2(k+1)}\sum_{j=1}^n\|u(t_j)\|_{k+2}^2\\\label{q11}
  &&+\tau\sum_{j=1}^n\|\bar{\partial}_tu(t_j)-u_t(t_j)\|_h^2
\end{eqnarray}
It is easy to see that
\begin{equation}\label{eq30}
\tau \sum_{j=1}^n\|u(t_j)\|_{k+2}^2\leq n\tau\max\limits_{1\leq j\leq n}\|u(t_j)\|_{k+2}^2\leq T\max\limits_{1\leq j\leq n}\|u(t_j)\|_{k+2}^2,
\end{equation}
and
\begin{equation}\label{q2}
\|\bar{\partial}_tu(t_j)-u_t(t_j)\|_h^2\leq C\tau\int_{t_{j-1}}^{t_j}\|u_{tt}(s)\|_h^2ds.
\end{equation}
Thus, substituting \eqref{eq30} and (\ref{q2}) into (\ref{q11})and using the discrete Gronwall lemma, we have
\begin{equation*}
    \|e_n^0\|_h^2\leq C \Big{(}h^{2(k+1)}\max\limits_{1\leq j\leq n}\|u(t_j)\|_{k+2}^2+{\tau}^2\int_0^{t_n}\|u_{tt}(s)\|_h^2ds\Big{)},
\end{equation*}
Then, by the triangle inequality, we complete the proof.\zb\par
By a similar argument to that of Theorem \ref{T:T2} and Theorem \ref{T:T3}, we can obtain the following error estimate for $\|u_x(t_n)-d_{w,r}U_n\|_h$. Due to space limitations, we omit the proof.
\begin{theorem}\label{T:T4}
Let $u$ and $U_n$ be the solution of problems (\ref{eq1}) and (\ref{eq15}), respectively, $u\in H^2(0,T;H^{1+r}(I))$ and $g_h=Q_hg$. Then, there exists a constant $C$ such that
\begin{eqnarray*}
\|u_x(t_n)-d_{w,r}U_n\|_h^2&\leq&C\Big{(}h^{2(k+1)}\max\limits_{1\leq j\leq n}\|u(t_j)\|_{k+2}^2\\
&+&h^{2(k+1)}\max\limits_{1\leq j\leq n}\|u_t(t_j)\|_{k+2}^2+\tau^2\int_0^{t_n}\|u_{tt}(s)\|_{k+2}^2ds\Big{)}.
\end{eqnarray*}
\end{theorem}
\section{Numerical experiment}
\begin{table}[h]
\begin{tabular}{ccccccccccc}
\hline
$x_j$&&\multicolumn{7}{c}{Numerical $U_n(x_j)$}&&Exact $u(x_j,t_n)$\\
\cline{3-9}
&&\multicolumn{3}{c}{$k=0$}&&\multicolumn{3}{c}{$k=1$}&&\\
\cline{3-5}
\cline{7-9}
&&\multicolumn{1}{c}{$N=80$}&&\multicolumn{1}{c}{$N=128$}&&\multicolumn{1}{c}{$N=80$}&&\multicolumn{1}{c}{$N=128$}&&\\
\hline
0.1&&0.20970&&0.21831&&0.22348&&0.22346&&0.22345\\
0.2&&0.42308&&0.43423&&0.43586&&0.43582&&0.43580\\
0.3&&0.61420&&0.62649&&0.62519&&0.62514&&0.62512\\
0.4&&0.76952&&0.78073&&0.77780&&0.77774&&0.77772\\
0.5&&0.87720&&0.87725&&0.87728&&0.87728&&0.87728\\
0.6&&0.90498&&0.90455&&0.90434&&0.90422&&0.90425\\
0.7&&0.84415&&0.83779&&0.83703&&0.83686&&0.83692\\
0.8&&0.67185&&0.65537&&0.65742&&0.65724&&0.65731\\
0.9&&0.38669&&0.35776&&0.36584&&0.36573&&0.36575\\
\hline
\end{tabular}
\caption{Comparison of the numerical solutions obtained with different weak finite element space $S_h$ and number of basis elements for $\nu=0.1$,
$\tau=0.0001$ at $t_n= 0.1$ with the exact solutions.}\label{Tab1}
\end{table}

In this section, we provide two numerical examples to confirm our theoretical analysis.
\subsection{Example 1}

\begin{table}[h]
\begin{tabular}{ccccccccccc}
\hline
$x_j$&&$t_n$&&\multicolumn{3}{c}{$\nu=0.1$}&&\multicolumn{3}{c}{$\nu=0.01$}\\
\cline{5-7}
\cline{9-11}
&&&&Numerical&&Exact&&Numerical&&Exact\\
&&&&$U_n(x_j)$&& $u(x_j,t_n)$&& $U_n(x_j)$&&$u(x_j,t_n)$\\
\hline
0.25&&0.4&&0.30892&&0.30889&&0.34194&&0.34191\\
    &&0.6&&0.24076&&0.24074&&0.26899&&0.26896\\
    &&0.8&&0.19569&&0.19568&&0.22150&&0.22148\\
    &&1.0&&0.16257&&0.16256&&0.18821&&0.18819\\
0.5&&0.4&&0.56966&&0.56963&&0.66074&&0.66071\\
   &&0.6&&0.44723&&0.44721&&0.52946&&0.52942\\
   &&0.8&&0.35926&&0.35924&&0.43917&&0.43914\\
   &&1.0&&0.29194&&0.29192&&0.37445&&0.37442\\
0.75&&0.4&&0.62542&&0.62544&&0.91021&&0.91026\\
    &&0.6&&0.48723&&0.48721&&0.76726&&0.76724\\
    &&0.8&&0.37394&&0.37392&&0.64743&&0.64740\\
    &&1.0&&0.28750&&0.28747&&0.55609&&0.55605\\
\hline
\end{tabular}
\caption{Comparison of the numerical solutions obtained for various values of $\nu$ and $N =80$, $\tau= 0.0001$ at
different times with the exact solutions.}\label{Tab2}
\end{table}
\begin{figure}
 \centering
 \includegraphics[width=0.5\textwidth]{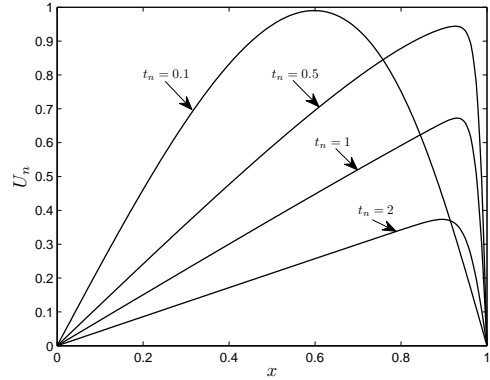}
 \caption{The numerical solutions by the WG finite element method with $k=1$, $N=80$, and $\tau=0.00001$ at different times for $\nu=0.01$.}\label{Fig1}
 \end{figure}
\begin{figure}
 \centering
 \includegraphics[width=0.5\textwidth]{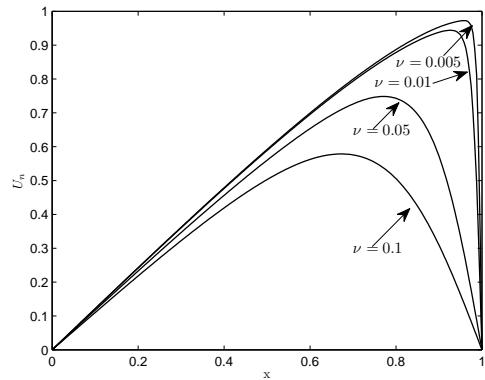}
 \caption{The numerical solutions by the WG finite element method with $k=1$, $N=80$, and $\tau=0.00001$ at $t_n=0.5$ for $\nu=0.1$, $\nu=0.05$, $\nu=0.01$, $\nu=0.005$.}\label{Fig2}
 \end{figure}
We first consider the problem defined by (\ref{eq1}) with initial condition
\begin{equation*}
g(x)=\sin(\pi x),\quad x\in [0,1].
\end{equation*}
The exact solution to this problem can be expressed as an infinite series
\begin{equation*}
    u(x,t)=2\pi \nu \frac{\sum\limits_{n=1}^{\infty}a_n\exp(-n^2\pi^2\nu t)n\sin(n\pi x)}{a_0+\sum\limits_{n=1}^{\infty}a_n\exp(-n^2\pi^2\nu t)\cos(n\pi)},
\end{equation*}
where the Fourier coefficients are defined by
\begin{equation*}
    a_0=\int_0^1\exp(-\frac{1-\cos(\pi x)}{2\pi\nu})dx,
\end{equation*}
and
\begin{equation*}
    a_n=2\int_0^1\exp(-\frac{1-\cos(\pi x)}{2\pi\nu})\cos(n\pi x)dx, \quad n=1,2,3,\cdots.
\end{equation*}\par
In order to implement the WG finite element schemes, we first divide $I=[0,1]$ into $N$ uniform intervals with length $h=\frac{1}{N}$. Then, we take $k=0$ and $k=1$ in the weak finite element space $S_h$ to carry out the WG finite element discretization (\ref{eq15}), respectively. \par
For the purpose of confirmation, the numerical solutions of the Burgers' equation obtained for various $\nu$ values at different times have been compared
with the exact solution. In Table \ref{Tab1}, a comparison of numerical and exact solutions for $k=0,1$, $N=80,128$, $\nu=0.1$,
$\tau=0.0001$ at $t_n = 0.1$ are shown. Table \ref{Tab2} shows a comparison of numerical and exact solutions for $k=1$, $N=80$, $\nu=0.1,0.01$,
$\tau=0.0001$ at $t_n = 0.4,0.6,0.8,1.0$. The profiles of the numerical solutions obtained by the WG finite element method with $k=1$, $N=80$, $\tau=0.00001$, $\nu=0.01$ at $t_n= 0.1,0.5,1.0,2.0$ are shown in Figure \ref{Fig1}. Figure \ref{Fig2} shows the profiles of the numerical solutions obtained by the WG finite element method with $k=1$, $N=80$, $\tau=0.00001$, $\nu=0.005,0.01,0.05,0.1$ at $t_n =0.5$.\par
As expected, the numerical solutions obtained by the WG finite element method provide high accuracy.

\subsection{Example 2}
\begin{figure}[ht]
 \centering
 \includegraphics[width=0.5\textwidth]{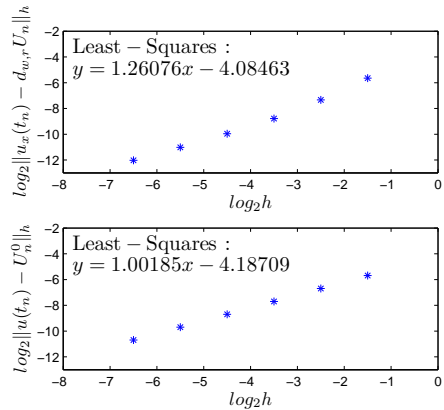}
 \caption{The convergence and error estimate of the WG finite element solution  with $k=0$, $\tau=0.00001$ at $t_n=1.0$ for $\nu=0.1$.}\label{Fig3}
 \end{figure}
 \begin{figure}[ht]
 \centering
 \includegraphics[width=0.5\textwidth]{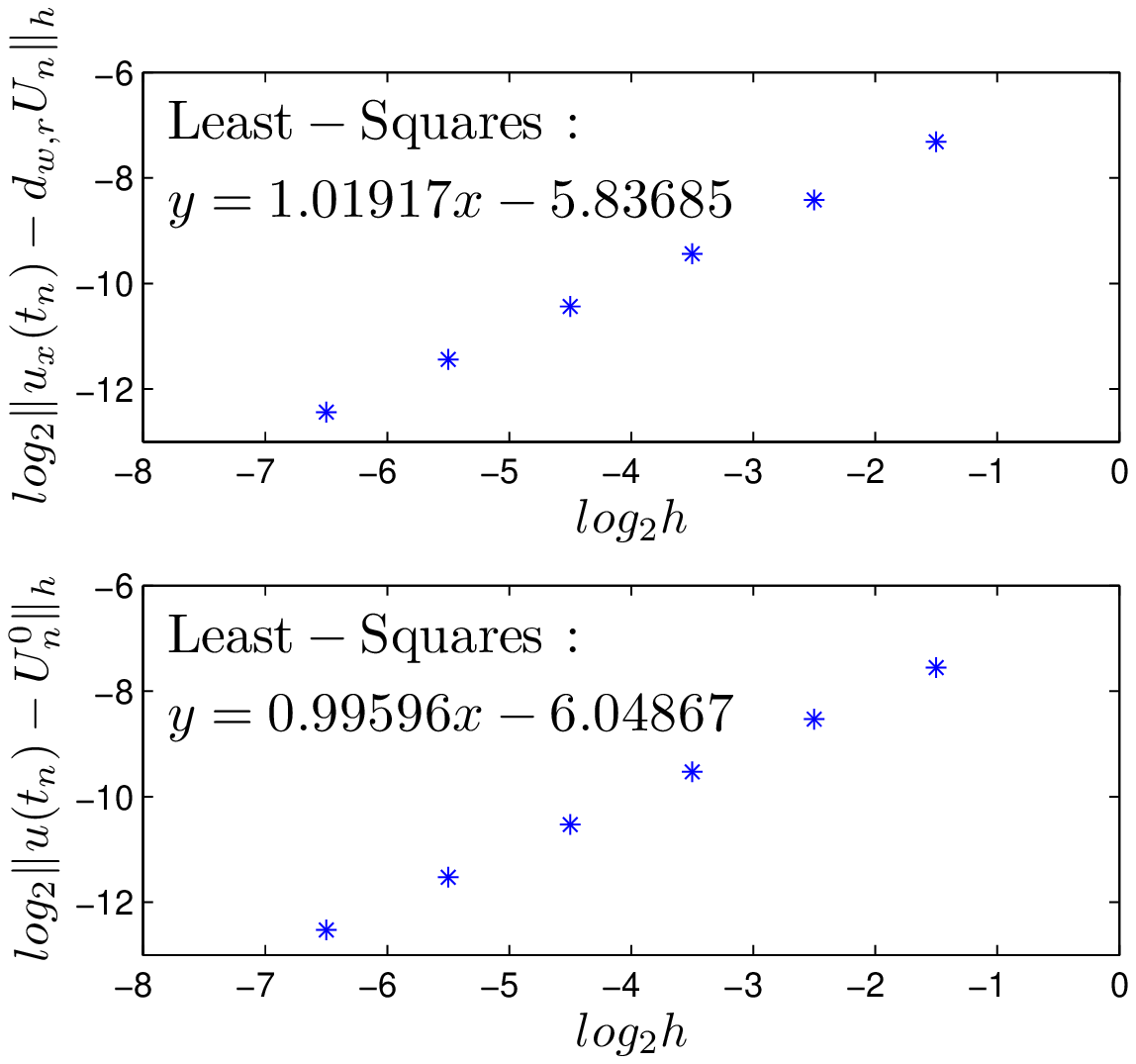}
 \caption{The convergence and error estimate of the WG finite element solution  with $k=0$, $\tau=0.00001$ at $t_n=1.0$ for $\nu=0.01$.}\label{Fig4}
 \end{figure}
  \begin{figure}[ht]
 \centering
 \includegraphics[width=0.5\textwidth]{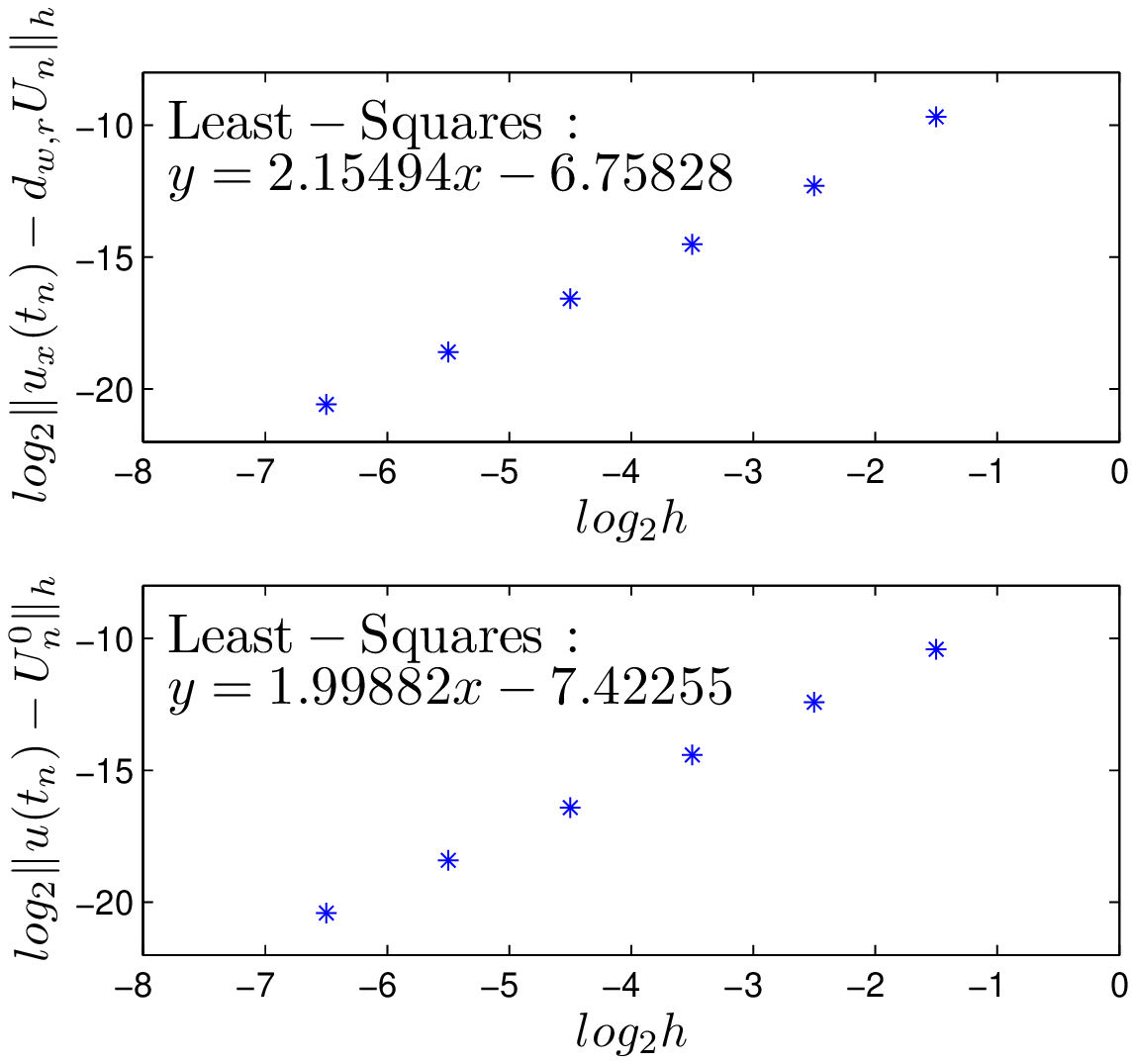}
 \caption{The convergence and error estimate of the WG finite element solution  with $k=1$, $\tau=0.00001$ at $t_n=1.0$ for $\nu=0.01$.}\label{Fig5}
 \end{figure}
   \begin{figure}[ht]
 \centering
 \includegraphics[width=0.5\textwidth]{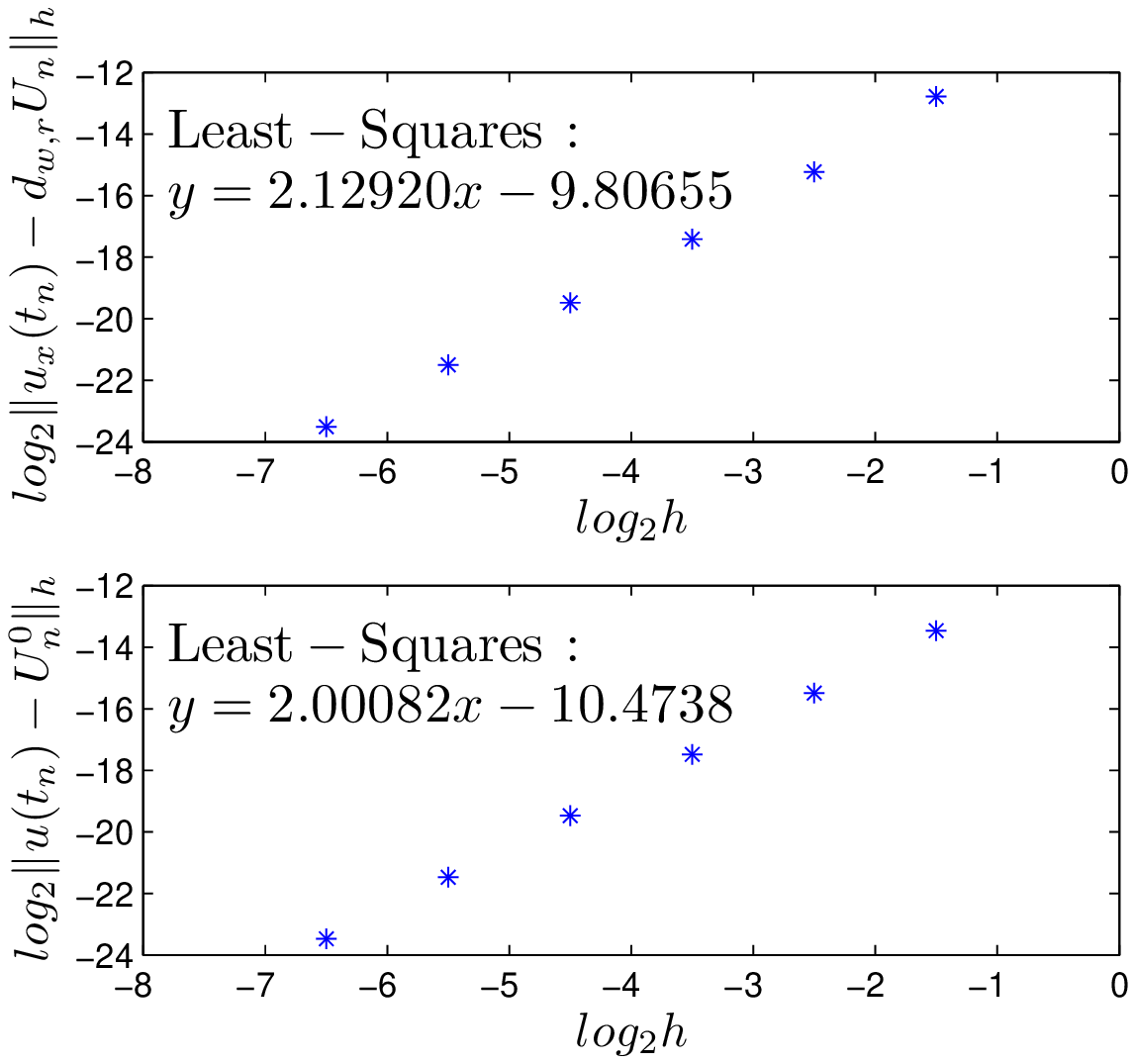}
 \caption{The convergence and error estimate of the WG finite element solution  with $k=1$, $\tau=0.00001$ at $t_n=1.0$ for $\nu=0.001$.}\label{Fig6}
 \end{figure}

In this example, we solve the Burgers' equation (\ref{eq1}) with initial condition
\begin{equation*}
    u(x,0)=2\nu\frac{\pi\sin(\pi x)}{\sigma+\cos(\pi x)},\quad x\in [0,1],
\end{equation*}
where $\sigma>1$ is a parameter. Here the initial condition is chosen such that the closed form of the exact solution is available, see \cite{b27},
\begin{equation*}
    u(x,t)=\frac{2\nu\pi e^{-\pi^2\nu t}\sin(\pi x)}{\sigma+e^{-\pi^2\nu t}\cos(\pi x)},\quad 0<x<1.
\end{equation*}
\par
Figure \ref{Fig3}-\ref{Fig6} show the convergence of the WG finite element method solutions obtained by $k=0,1$ and $\tau=0.00001$ at $t_n=1.0$ for $\nu=0.1,0.01,0.001$. Clearly, it is shown that the WG finite element solution exhibits $k+1$ order convergence in discrete $H^1$-norm and $L^2$-norm, which are consistent with Theorem \ref{T:T3} and Theorem \ref{T:T4}.

\end{document}